# Potentials for a multidimensional elliptic equation with one line of degeneration and their applications to boundary value problems

H. M. Srivastava, A. Hasanov, T. G. Ergashev

H. M. Srivastava, Department of Mathematics and Statistics, University of Victoria, Victoria, British Columbia V8W 3R4, Canada *and* Department of Medical Research, China Medical University Hospital, China Medical University, Taichung 40402, Taiwan, Republic of China. E-Mail: harimsri@math.uvic.ca
Anvar Hasanov, Institute of Mathematics. 81 Mirzo-Ulugbek Street, Tashkent 700170, Uzbekistan. E-Mail: anvarhasanov@yahoo.com
T. G. Ergashev, Institute of Mathematics. 81 Mirzo-Ulugbek Street, Tashkent 700170, Uzbekistan. E-Mail: ergashev.tukhtasin@gmail.com

**Abstract:**
Potentials play an important role in solving boundary value problems for elliptic equations. In the middle of the last century, a potential theory was constructed for a two-dimensional elliptic equation with one singular coefficient. In the study of potentials, the properties of the fundamental solutions of the given equation are essentially used. At the present time, fundamental solutions of a multidimensional elliptic equation with one degeneration line are already known. In this paper, we investigate the potentials of the double- and simple-layers for this equation, with the help of which limit theorems are proved and integral equations containing in the kernel the density of the above potentials are derived.

## 1. Introduction

Potential theory has played a paramount role in both analysis and computation for boundary value problems for elliptic equations. Numerous applications can be found in mechanics, fluid mechanics, elastodynamics, electromagnetic, and acoustics. Results from potential theory allow us to represent boundary value problems in integral equation form. For problems with known Green's functions, an integral equation formulation leads to powerful numerical approximation schemes.

The double-layer and simple-layer potentials play an important role in solving boundary value problems for elliptic equations. For example, the representation of the solution of the (first) boundary value problem is sought as a double-layer potential with unknown density and an application of certain property leads to a Fredholm equation of the second kind for determining the function (see [10] and [14]).

Let be $R_m^+$ – the half-space $x_1 > 0$ of the $m-$ dimensional Euclidean space of points $x \equiv (x_1,...,x_m)$.

We consider the equation

$$H_\alpha^m(u) \equiv \sum_{i=1}^m u_{x_i x_i} + \frac{2\alpha}{x_1} u_{x_1} = 0, \qquad (1.1)$$

where $\alpha$ and $m$ are constants, and $0 < 2\alpha < 1$, $m \geq 2$.

Various interesting problems associated with the equation $H_\alpha^2(u) = 0$ were studied by many authors (see [1-3], [5-7], [11, 12], [16, 17], [19]). For example, by applying a method of complex analyses (based upon analytic functions), Gilbert [8] constructed an integral representation of solutions of generalized Helmholtz equation. Weinstein first considered fractional dimensional space in potential theory [20, 21]. Relatively few papers have been devoted to the potential theory for equation (1.1) at $m > 2$. We note only the papers [13, 15].

Fundamental solutions of the equation (1.1) were constructed recently (see [18]) and they are expressed by hypergeometric function of Gauss

$$F(a,b;c;z) = \sum_{n=0}^{\infty} \frac{(a)_n (b)_n}{n!(c)_n} z^n,$$

where $(\kappa)_\nu$ denotes the general Pochhammer symbol or the shifted factorial, since $(1)_n = n!$ $(n \in N_0 := N \cup \{0\}; N := \{1,2,3,...\})$, which is defined (for $\kappa, \nu \in C$), in terms of the familiar Gamma function, by

$$(\kappa)_\nu := \frac{\Gamma(\kappa+\nu)}{\Gamma(\kappa)} = \begin{cases} 1 & (\nu = 0; \kappa \in C \setminus \{0\}) \\ \kappa(\kappa+1)...(\kappa+n-1) & (\nu = n \in N; \kappa \in C), \end{cases}$$

it being understood conventionally that $(0)_0 := 1$ and assumed tacitly that the $\Gamma-$ quotient exists.

The fundamental solutions of equation (1.1) look like [18]:

$$q_1(x,\xi) = k_1 (r^2)^{-\alpha - \frac{m-2}{2}} F\left(\alpha + \frac{m-2}{2}, \alpha; 2\alpha; \zeta\right), \tag{1.2}$$

$$q_2(x,\xi) = k_2 (r^2)^{\alpha - \frac{m}{2}} x_1^{1-2\alpha} \xi_1^{1-2\alpha} F\left(\frac{m}{2} - \alpha, 1-\alpha; 2-2\alpha; \zeta\right), \tag{1.3}$$

where

$$\xi = (\xi_1, \xi_2, ..., \xi_m), \quad r^2 = \sum_{i=1}^{m}(x_i - \xi_i)^2, \quad r_1^2 = (x_1 + \xi_1)^2 + \sum_{i=2}^{m}(x_i - \xi_i)^2, \tag{1.4}$$

$$\zeta = 1 - \frac{r_1^2}{r^2}, \quad k_1 = \frac{\Gamma(\alpha)\Gamma(\alpha + (m-2)/2)}{4^{1-\alpha} \pi^{m/2} \Gamma(2\alpha)}, \quad k_2 = \frac{\Gamma(1-\alpha)\Gamma(-\alpha + m/2)}{4^\alpha \pi^{m/2} \Gamma(2-2\alpha)}.$$

The fundamental solutions given by (1.2) and (1.3) posses the following properties:

$$\left.\frac{\partial q_1(x,\xi)}{\partial x_1}\right|_{x_1=0} = 0, \tag{1.5}$$

$$q_2(x,\xi)\big|_{x_1=0} = 0.$$

Here, by making use of the fundamental solution given by (1.2) in the domain $\Omega$ defined by

$$\Omega \subset R_m^+ := \{x : x_1 > 0\}, \tag{1.6}$$

we aim at investigating a double-layer and simple-layer potentials for the equation (1.1). Furthermore, we prove some results on limiting values of this potentials. These results are used in the last section of this paper, where Holmgren's problem is investigated in more general domains. Throughout this paper it is assumed that the dimension of the space $m > 2$.

## 2. Green's Formula

We consider the identity

$$x_1^{2\alpha} \left[ u H_\alpha^m(v) - v H_\alpha^m(u) \right] = \sum_{i=1}^{m} \frac{\partial}{\partial x_i} \left[ x_1^{2\alpha} \left( v_{x_i} u - v u_{x_i} \right) \right]. \tag{2.1}$$

Integrating both sides of the last identity in a domain $\Omega$ located and bounded in the half-space ($x_1 > 0$), and using the Gauss-Ostrogradsky formula, we obtain

$$\int_\Omega x_1^{2\alpha} \left[ u H_\alpha^m(v) - v H_\alpha^m(u) \right] dx = \int_S x_1^{2\alpha} \sum_{i=1}^{m} \left( u \frac{\partial v}{\partial x_i} - v \frac{\partial u}{\partial x_i} \right) \cos(n, x_i) dS, \tag{2.2}$$

where $S$ is a boundary of $\Omega$.

The Green formula (2.2) is derived under the following assumptions: the functions $u(x)$, $v(x)$, and their first-order partial derivatives are continuous in the closed domain $\bar\Omega$, the second-order partial derivatives are continuous inside $\Omega$ and the integrals by $\Omega$ which contain $H_\alpha^m(u)$ and $H_\alpha^m(v)$ have a meaning. If $H_\alpha^m(u)$ and $H_\alpha^m(v)$ do not have continuity up to $S$, then they are improper integrals that are obtained as limits on any sequence of domains $\Omega_k$ that are contained inside $\Omega$ when these domains $\Omega_k$ tend to $\Omega$, so that any point inside $\Omega_k$ is inside $\Omega$ the domains, starting with some number.

If $u$ and $v$ are solutions of equation (1.1), then from formula (2.2) we have

$$\int_S x_1^{2\alpha}\left(u\frac{\partial v}{\partial n} - v\frac{\partial u}{\partial n}\right)dS = 0. \tag{2.3}$$

Here

$$\frac{\partial}{\partial n} = \sum_{i=1}^m \cos(n, x_i) \cdot \frac{\partial}{\partial x_i} \tag{2.4}$$

is a normal derivative, and $n$ is the outer normal to the surface $S$.

Assuming $v = 1$ in (2.2) and replacing $u$ by $u^2$, we obtain

$$\int_\Omega x_1^{2\alpha} \sum_{i=1}^m \left(\frac{\partial u}{\partial x_i}\right)^2 dx = \int_S x_1^{2\alpha} u \frac{\partial u}{\partial n} dS, \tag{2.5}$$

where $u(x)$ is the solution of equation (1.1)

The special case of (2.3) when $v = 1$ reduces to the following form:

$$\int_S x_1^{2\alpha} \frac{\partial u}{\partial n} dS = 0. \tag{2.6}$$

We note from (2.6) that the integral of the normal derivative of the solution of equation (1.1) with a weight $x_1^{2\alpha}$ along the boundary $S$ of the domain $\Omega$ in (1.6) is equal to zero.

## 3. A double-Layer Potential $w^{(1)}(x)$.

The surface $\Gamma$ in the Euclidean space $E_m$, that satisfies the following two conditions is called the *Lyapunov surface*:

(i) At any point of the surface $\Gamma$ there is a definite normal.

(ii) Let $x$ and $\xi$ be points of the surface $\Gamma$, and $\vartheta$ angle between these normals. There exist positive constants $a$ and $\alpha$, such that

$\vartheta \leq ar^\alpha$.

(iii) With respect to the surface $\Gamma$ we shall assume that it approaches the hyperplane $x_1 = 0$ under right angle.

Let $\Omega$ be a finite domain in $R_m^+$, bounded by the open part $\Gamma_1$ of the hyper plane $x_1 = 0$ and the Lyapunov surface $\Gamma$. The boundary of the domain $\Gamma_1$ will be denoted by $\gamma$.

We consider the following integral

$$w^{(1)}(x) = \int_\Gamma \xi_1^{2\alpha} \mu_1(\xi) \frac{\partial q_1(\xi, x)}{\partial n_\xi} d_\xi \Gamma, \tag{3.1}$$

where the denseness $\mu_1(x) \in C(\overline{\Gamma})$ and $q_1(\xi, x)$ is given in (1.2). We call the integral (3.1) *a double-layer potential with denseness* $\mu_1(\xi)$. When $\mu_1(\xi) = 1$, we denote the double-layer potential (3.1) by $w_1^{(1)}(x)$.

We now investigate some properties of a double-layer potential $w_1^{(1)}(x)$.

**Lemma 1.** *The following formula holds true:*

$$w_1^{(1)}(x) = \begin{cases} -1, & x \in \Omega, \\ -\dfrac{1}{2}, & x \in \Gamma, \\ 0, & x \notin \bar{\Omega}, \end{cases}$$

*where a domain $\Omega$ and the surface $\Gamma$ are described as in this section and $\bar{\Omega} := \Omega \cup \Gamma$.*

**Proof. Case 1.** When $x \in \Omega$, we cut a ball centered at $x$ with a radius $\rho$ off the domain $\Omega$ and denote the remaining part by $\Omega_\rho$ and the sphere of the cut-off-ball by $C_\rho$. The function $q_1(\xi; x)$ in (1.2) is a regular solution of the equation (1.1) in the domain $\Omega_\rho$.

Using the well-known derivative formula of Gauss hypergeometric function [4], we have

$$\frac{\partial q_1(\xi, x)}{\partial \xi_1} = -k_1(2\alpha + m - 2)(\xi_1 - x_1)(r^2)^{-\alpha - \frac{m}{2}} F\left(\alpha + \frac{m}{2}, \alpha; 2\alpha; \zeta\right)$$
$$-(2\alpha + m - 2)k_1 x_1 (r^2)^{-\alpha - \frac{m}{2}} F\left(\alpha + \frac{m}{2}, 1 + \alpha; 1 + 2\alpha; \zeta\right). \tag{3.2}$$

$$\frac{\partial q_1(\xi, x)}{\partial \xi_i} = -k_1(2\alpha + m - 2)(\xi_i - x_i)(r^2)^{-\alpha - \frac{m}{2}} F\left(\alpha + \frac{m}{2}, \alpha; 2\alpha; \zeta\right), \quad i = 2, \ldots, m. \tag{3.3}$$

Using (3.2) and (3.3), by virtue of (2.4), we find

$$\frac{\partial}{\partial n_\xi}\{q_1(\xi, x)\} = (2\alpha + m - 2)k_1 (r^2)^{-\alpha - \frac{m-3}{2}} F\left(\alpha + \frac{m}{2}, \alpha; 2\alpha; \zeta\right) \frac{\partial}{\partial n_\xi}\left(\frac{1}{r}\right)$$
$$-(2\alpha + m - 2)k_1 x_1 (r^2)^{-\alpha - \frac{m}{2}} F\left(\alpha + \frac{m}{2}, 1 + \alpha; 1 + 2\alpha; \zeta\right) \cos(n, \xi_1). \tag{3.4}$$

Applying (2.6) and in view of (3.1), we get following formula

$$w_1^{(1)}(x) = \lim_{\rho \to 0} \int_{C_\rho} \xi_1^{2\alpha} \frac{\partial}{\partial n_\xi}\{q_1(\xi, x)\} d_\xi C_\rho - \int_{\Gamma_1} \left[\xi_1^{2\alpha} \frac{\partial}{\partial n_\xi}\{q_1(\xi, x)\}\right]_{\xi_1 = 0} d_{\xi'} \Gamma_1. \tag{3.5}$$

Substituting from (3.4) into (3.5), by virtue of (1.5), we find that
$$w_1^{(1)}(x) = i_1(x) + j_1(x), \tag{3.6}$$
where

$$i_1(x) = (2\beta + m - 2)k_1 \int_{C_\rho} \xi_1^{2\alpha} (r^2)^{-\alpha - \frac{m-3}{2}} F\left(\alpha + \frac{m}{2}, \alpha; 2\alpha; \zeta\right) \frac{\partial}{\partial n_\xi}\left(\frac{1}{r}\right) d_\xi C_\rho,$$

$$j_1(x) = -(2\beta + m - 2)k_1 x_1 \int_{C_\rho} \xi_1^{2\alpha} (r^2)^{-\alpha - \frac{m}{2}} F\left(\alpha + \frac{m}{2}, 1 + \alpha; 1 + 2\alpha; \zeta\right) \cos(n, \xi_1) d_\xi C_\rho.$$

Now, introducing the spherical coordinates:
$$\xi_1 = x_1 + \rho \cos \varphi_1,$$
$$\xi_i = x_i + \rho \sin \varphi_1 \sin \varphi_2 \ldots \sin \varphi_{i-1} \cos \varphi_i, i = 2, \ldots, m-1, \tag{3.7}$$
$$\xi_m = x_m + \rho \sin \varphi_1 \sin \varphi_2 \ldots \sin \varphi_{m-2} \sin \varphi_{m-1},$$
$$(0 \leq \rho \leq r, 0 \leq \varphi_1 \leq \pi, \ldots, 0 \leq \varphi_{m-2} \leq \pi, 0 \leq \varphi_{m-1} \leq 2\pi)$$
we get

$$i_1(x) = -(2\beta + m - 2)k_1 \rho^{-2\beta} \int_0^{2\pi} d\varphi_{m-1} \int_0^{\pi} \sin\varphi_{m-2} d\varphi_{m-2} \int_0^{\pi} \sin^2\varphi_{m-3} d\varphi_{m-3} \ldots \int_0^{\pi} \sin^{m-2}\varphi_1$$
$$\times (x_m + \rho\sin\varphi_1 \sin\varphi_2 \ldots \sin\varphi_{m-2}\sin\varphi_{m-1})^{2\beta} F\left(\beta + \frac{m}{2}, \beta; 2\beta; \zeta_\rho\right) d\varphi_1, \quad (3.8)$$

where
$$\zeta_\rho = -4x_m \rho^{-2}(x_m + \rho\sin\varphi_1 \sin\varphi_2 \ldots \sin\varphi_{m-2}\sin\varphi_{m-1}).$$

Let us investigate the integrand in (3.8). Applying the well-known formula

$$F(a,b;c;x) = (1-x)^{-b} F\left(c-a,b;c;\frac{x}{x-1}\right), \quad (3.9)$$

we obtain

$$i_1(x) = -(2\beta + m - 2)k_1 \int_0^{2\pi} d\varphi_{m-1} \int_0^{\pi} \sin\varphi_{m-2} d\varphi_{m-2} \int_0^{\pi} \sin^2\varphi_{m-3} d\varphi_{m-3} \ldots \int_0^{\pi} \sin^{m-2}\varphi_1$$

$$\times \frac{(x_m + \rho\sin\varphi_1 \sin\varphi_2 \ldots \sin\varphi_{m-2}\sin\varphi_{m-1})^{2\beta}}{\left[\rho^2 + 4x_m(x_m + \rho\sin\varphi_1 \sin\varphi_2 \ldots \sin\varphi_{m-2}\sin\varphi_{m-1})\right]^\beta}$$

$$\times F\left(\beta - \frac{m}{2}, \beta; 2\beta; \frac{4x_m(x_m + \rho\sin\varphi_1 \sin\varphi_2 \ldots \sin\varphi_{m-2}\sin\varphi_{m-1})}{\rho^2 + 4x_m(x_m + \rho\sin\varphi_1 \sin\varphi_2 \ldots \sin\varphi_{m-2}\sin\varphi_{m-1})}\right) d\varphi_1.$$

Using the well-known Gauss's summation formula for $F$ (see [4], p.104, eq.(46))
$$F = (a,b;c;1) = \frac{\Gamma(c)\Gamma(c-a-b)}{\Gamma(c-a)\Gamma(c-b)}, \quad \text{Re}(c-a-b) > 0; \ c \neq 0, -1, -2, \ldots$$

we obtain

$$\lim_{\rho \to 0} i_1(x) = -2^{1-2\alpha} k_1 L_m \frac{\Gamma(2\alpha)\Gamma\left(\frac{m}{2}\right)}{\Gamma\left(\alpha - 1 + \frac{m}{2}\right)\Gamma(\alpha)},$$

where

$$L_m = \int_0^{2\pi} d\varphi_{m-1} \int_0^{\pi} \sin\varphi_{m-2} d\varphi_{m-2} \int_0^{\pi} \sin^2\varphi_{m-3} d\varphi_{m-3} \ldots \int_0^{\pi} \sin^{m-2}\varphi_1 d\varphi_1.$$

With the help of elementary transformations, known from the course of mathematical analysis, it is not difficult to establish that

$$L_{2m} = \frac{2\pi^m}{(m-1)!}, \qquad L_{2m+1} = \frac{2^{m+1}\pi^m}{(2m-1)!!}, \quad m = 1, 2, 3, \ldots$$

Now, using the Legendre's duplication formula (see [4], p.5, formula (15))
$\Gamma(2a) = 2^{2a-1}\pi^{-1/2}\Gamma(a)\Gamma(a+1/2)$, we have

$$\lim_{\rho \to 0} i_1(x) = -1. \quad (3.10)$$

Similarly, by considering the corresponding identities we find that
$$\lim_{\rho \to 0} j_1(x) = 0. \quad (3.11)$$

Hence, by view of (3.10) and (3.11), the formula (3.6) in the case of $x \in \Omega$ becomes
$w_1^{(1)}(x) = -1.$

**Case 2.** When $x \in \Gamma$, we cut we cut a sphere $C_\rho$ centered at $x$ with a radius $\rho$ off the domain $\Omega$ and denote the remaining part of the surface by $\Gamma'$, that is $\Gamma' = \Gamma - \Gamma_\rho$. Let $C'_\rho$ denote a part of the sphere $C_\rho$ lying inside the domain $\Omega$. We consider the domain $\Omega_\rho$ which is bounded by a surface $\Gamma'$, $C'_\rho$ and $\Gamma_1$ and its boundary $\gamma$. Then we have

$$w_1^{(1)}(x) = \lim_{\rho \to 0} \int_{\Gamma'} \xi_1^{2\alpha} \frac{\partial q_1(\xi, x)}{\partial n_\xi} d_\xi \Gamma'. \tag{3.12}$$

When the point $x$ lies outside the domain $\Omega_\rho$, it is found that, in this domain $q_1(\xi, x)$ is a regular solution of the equation (1.1). Therefore, by virtue of (2.6), we have

$$\int_{\Gamma'} \xi_1^{2\alpha} \frac{\partial q_1(\xi, x)}{\partial n_\xi} d_\xi \Gamma' = \int_{C'_\rho} \xi_1^{2\alpha} \frac{\partial q_1(\xi, x)}{\partial n_\xi} d_\xi C'_\rho. \tag{3.13}$$

Substituting from (3.13) into (3.12), we get

$$w_1^{(1)}(x) = \lim_{\rho \to 0} \int_{C'_\rho} \xi_1^{2\alpha} \frac{\partial q_1(\xi, x)}{\partial n_\xi} d_\xi C'_\rho.$$

Similarly, by again introducing the spherical coordinates (3.7) centered at the point $x$, we have

$$w_1^{(1)}(x) = -\frac{1}{2}, \ x \in \Gamma.$$

**Case 3.** When $x \notin \overline{\Omega}$, it is noted that the function $q_1(\xi, x)$ is a regular solution of the equation (1.1). Hence, in view of formula (2.6), we have $w_1^{(1)}(x) = 0, \ x \notin \overline{\Omega}$.

The proof of Lemma 1 is thus completed.

**Lemma 2.** *The following formula holds true:*

$$w_1^{(1)}(0, x') = \begin{cases} -1, & (0, x') \in \Gamma_1, \\ -\dfrac{1}{2}, & (0, x') \in \gamma, \\ 0, & (0, x') \notin \Gamma_1 \cup \gamma, \end{cases}$$

where $x' = (x_2, x_2, ..., x_m)$.

**Proof.** For considering the first case when $(0, x') \in \Gamma_1$, we introduce a hyper plane $x_1 = \delta$ for a sufficiently small positive real number $\delta$ and consider a domain $\Omega_\delta$ which is the part of the domain $\Omega$ lying above the hyper plane $x_1 = \delta$. Applying the formula (2.6), we obtain

$$w_1^{(1)}(0, x') = \int_{\Gamma_{1\delta}} \delta^{2\alpha} \frac{\partial q_1(\xi, 0, x')}{\partial \xi_1}\bigg|_{\xi_1 = \delta} d\xi' + \int_{\Gamma_\delta} \xi_1^{2\alpha} \frac{\partial q_1(\xi, 0, x')}{\partial n_\xi} d_\xi \Gamma_\delta, \tag{3.14}$$

where $\xi' = (\xi_2, \xi_2, ... \xi_m)$ and $\Gamma_\delta$ is a part of the surface $\Gamma$ at $\xi_1 \leq \delta$. It follows from (3.2) and (3.14) that

$$w_1^{(1)}(0, x') = -k_1(2\alpha + m - 2) \lim_{\delta \to 0} \int_{\Gamma_{1\delta}} \delta^{1+2\alpha} \left( \sum_{i=2}^{m} (\xi_i - x_i)^2 + \delta^2 \right)^{-\alpha - \frac{m}{2}} d\xi'. \tag{3.15}$$

Now, we transform the expression (3.15). Instead of $\xi'$, we introduce new integration variables $t' = \dfrac{\xi' - x'}{\delta}$, where $t' = (t_2, t_2, ..., t_m)$. Making the change of variables, after passing to limit, we obtain

$$w_1^{(1)}(0, x') = -k_1(2\alpha + m - 2) \int_{-\infty}^{\infty} ... \int_{-\infty}^{\infty} \left( 1 + \sum_{i=2}^{m} t_i^2 \right)^{-\alpha - \frac{m}{2}} dt'.$$

It is known that [9, p.637, formula 4.638.3]

$$\int_{-\infty}^{\infty}...\int_{-\infty}^{\infty}\left(1+\sum_{i=2}^{m}t_i^2\right)^{-\alpha-\frac{m}{2}}dt' = \pi^{(m-1)/2}\Gamma\left(\alpha+\frac{1}{2}\right)\Gamma^{-1}\left(\alpha+\frac{m}{2}\right).$$

Using this formula, we have
$w_1^{(1)}(0,x') = -1$, $(0,x') \in \Gamma_1$.

The other cases can be proved by using arguments similar to those detailed above in the first case. This evidently completes our proof of Lemma 2.

**Theorem 1.** *For any points $x$ and $\xi \in R_m^+$ and $\xi \neq x$, the following inequality holds true:*

$$|q_1(\xi,x)| \leq k_1 \frac{r_1^{-2\alpha}}{r^{m-2}}, \tag{3.16}$$

*where $m > 2$,, $\alpha$ is real parameter with $0 < 2\alpha < 1$ as in the equation (1.1), $r$ and $r_1$ are as in (1.4).*

**Proof.** Using the formula (3.11), we obtain

$$q_1(x,\xi) = k_1 (r^2)^{-\frac{m-2}{2}}(r_1^2)^{-\alpha} F\left(\alpha - \frac{m-2}{2}, \alpha; 2\alpha; 1 - \frac{r^2}{r_1^2}\right).$$

*This immediately implies the estimate (3.16). Theorem 1 is proved.*

**Theorem 2.** *If a surface $\Gamma$ satisfies the conditions (i)-(iii), then the following inequality holds true*

$$\int_\Gamma \xi_1^{2\alpha} \left|\frac{\partial q_1(\xi,x)}{\partial n_\xi}\right| d_\xi \Gamma \leq C,$$

*where $C$ is a constant.*

**Proof.** Theorem 2 follows from conditions (i)-(iii) and (3.4).

**Theorem 3.** *The following limiting formulas hold true for a double-layer potential (3.1):*

$$w_i^{(1)}(t) = -\frac{1}{2}\mu_1(t) + \int_\Gamma \mu_1(s) K_1(s,t) d_s \Gamma, \tag{3.17}$$

and

$$w_e^{(1)}(t) = \frac{1}{2}\mu_1(t) + \int_\Gamma \mu_1(s) K_1(s,t) d_s \Gamma, \tag{3.18}$$

*where* $\mu_1(t) \in C(\overline{\Gamma})$, $s = (s_1, s_2, ..., s_m)$, $t = (t_1, t_2, ..., t_m)$,

$$K_1(s,t) = s_1^{2\alpha} \frac{\partial}{\partial n_s}\{q_1(s,t)\},\ s \in \Gamma,\ t \in \Gamma, \tag{3.19}$$

$w_i^{(1)}(t)$ *and* $w_e^{(1)}(t)$ *are limiting values of the double-layer potential (3.1) at $t \to \Gamma$ from the inside and the outside, respectively.*

**Proof.** We find from Lemma 1, in conjunction with Theorems 1 and 2, that each of the limiting formulas asserted by Theorem 3 holds true.

**4.  Simple –Layer Potential $v^{(1)}(x)$**

We consider the following integral
$$v^{(1)}(x) = \int_\Gamma \rho_1(\xi) q_1(\xi,x) d_\xi \Gamma, \tag{4.1}$$

where the density $\rho_1(x) \in \overline{\Gamma}$ and $q_1(\xi, x)$ is given in (1.2). We call the integral (4.1) *a simple-layer potential with denseness* $\rho_1(\xi)$. When $\rho_1(\xi) = 1$, we denote the simple-layer potential (4.1) by $v_1^{(1)}(x)$.

We now investigate some properties of a simple-layer potential $v_1^{(1)}(x)$.

The simple-layer potential (4.1) is defined throughout the half-space $x_1 > 0$ and is a continuous function when passing through the surface $\Gamma$. Obviously, a simple-layer potential is a regular solution of equation (1.1) in any domain lying in the half-space $x_1 > 0$. It is easy that as the point $x$ tends to infinity, a simple layer potential $v^{(1)}(x)$ tends to zero. Indeed, let the point $x$ be on the hemisphere

$$\sum_{i=1}^{m} x_i^2 = R^2, x_1 > 0,$$

then by virtue of (1.2) and (1.3) we have

$$\left|v^{(1)}(x)\right| \leq \int_\Gamma |\rho_1(t)||q_1(t,x)|d_t\Gamma \leq MR^{2-2\alpha-m} \quad (R \geq R_0)$$

where $M$ is a constant.

We take an arbitrary point $N(x_0)$ on the surface $\Gamma$ and draw a normal at this point. Consider on this normal some point $M(x)$, not lying on the surface $\Gamma$, we find the normal derivative of the simple-layer potential (4.1):

$$\frac{\partial}{\partial n}\{v^{(1)}(x)\} = \int_\Gamma \rho_1(t) \frac{\partial}{\partial n_x}\{q_1(t,x)\} d_t\Gamma.$$

The integral (4.4) exists also in the case when the point $M(x)$ coincides with the point $N$, mentioned above.

**Theorem 3.** *The following limiting formulas hold true for a simple-layer potential (4.1):*

$$t_1^{2\alpha} \frac{\partial}{\partial n}\{v^{(1)}(t)\}_i = \frac{1}{2}\rho_1(t) + \int_\Gamma \rho_1(s) K_1(s,t) d_s\Gamma, \quad (4.2)$$

*and*

$$t_1^{2\alpha} \frac{\partial}{\partial n}\{v^{(1)}(t)\}_e = -\frac{1}{2}\rho_1(t) + \int_\Gamma \rho_1(s) K_1(s,t) d_s\Gamma, \quad (4.3)$$

*where* $\rho_1(t) \in C(\overline{\Gamma})$,

$K_1(s,t) = t_1^{2\alpha} \dfrac{\partial}{\partial n_t}\{q_1(s,t)\}, s \in \Gamma, t \in \Gamma,$

$t_1^{2\alpha} \dfrac{\partial}{\partial n}\{v^{(1)}(t)\}_i$ *and* $t_1^{2\alpha} \dfrac{\partial}{\partial n}\{v^{(1)}(t)\}_e$ *are limiting values of the normal derivative of simple-layer potential (4.1) at* $t \to \Gamma$ *from the inside and the outside, respectively.*

From these formulas the jump in the normal derivative of the simple-layer potential follows immediately:

$$\frac{\partial}{\partial n}\{v^{(1)}(t)\}_i - \frac{\partial}{\partial n}\{v^{(1)}(t)\}_e = t_1^{-2\alpha}\rho_1(t). \quad (4.4)$$

For the future it is useful to note that when the point tends to infinity,

$$\left|\frac{\partial}{\partial n_s}\{v^{(1)}(s)\}_i\right| \le NR^{1-2\alpha-m}, \quad (R \ge R_0)$$

where $N$ is a constant.

In exactly the same way as in the derivation of (2.5), it is not difficult to show that Green's formulas are applicable to the simple-layer potential (4.1):

$$\int_\Omega x_1^{2\alpha} \sum_{k=1}^m \left(\frac{\partial v^{(1)}}{\partial x_k}\right) dx = \int_S x_1^{2\alpha} v^{(1)}(x) \frac{\partial}{\partial n}\{v^{(1)}(x)\}_i \, dS, \qquad (4.5)$$

$$\int_\Omega x_1^{2\alpha} \sum_{k=1}^m \left(\frac{\partial v^{(1)}}{\partial x_k}\right) dx = -\int_S x_1^{2\alpha} v^{(1)}(x) \frac{\partial}{\partial n}\{v^{(1)}(x)\}_e \, dS. \qquad (4.6)$$

## 5. Integral Equations For Denseness's

Formulas (3.17), (3.18), (4.2) and (4.5) can be written as integral equations for denseness's:

$$\mu_1(s) - \lambda \int_\Gamma K_1(s,t) \mu_1(t) d_t\Gamma = f_1(s), \qquad (5.1)$$

$$\rho_1(s) - \lambda \int_\Gamma K_1(t,s) \rho_1(t) d_t\Gamma = g_1(s), \qquad (5.2)$$

where

$$\lambda = 2, \; f_1(s) = -2w_i^{(1)}(s), \; g_1(s) = -2\frac{\partial}{\partial n}\{v^{(1)}(s)\}_e,$$

$$\lambda = -2, \; f_1(s) = 2w_e^{(1)}(s), \; g_1(s) = 2\frac{\partial}{\partial n}\{v^{(1)}(s)\}_i.$$

Equations (5.1) and (5.2) are conjugated and, by Theorem 1, Fredholm theory is applicable to them. We show that $\lambda = 2$ is not an eigenvalue of the kernel $K_1(s,t)$. This assertion is equivalent to the fact that the homogeneous integral equation

$$\rho_1(t) - 2\int_\Gamma K_1(s,t) \rho_1(s) d_s\Gamma = 0 \qquad (5.3)$$

has no nontrivial solutions.

Let $\rho_0(t)$ be a continuous non-trivial solution of equation (5.3). We note that all bounded solutions of equation (5.3) satisfy condition $|\rho(t)| \le Mt_1^{2\alpha}$. The simple-layer potential with denseness $\rho_0(t)$ gives us a function $v_0(x)$ which is a solution of equation (1.1) in domains $\Omega$ and $R_m^+ \setminus \overline{\Omega}$. The limiting values of the normal derivative of $\frac{\partial}{\partial n}\{v_0(s)\}_e$ are zero, by virtue of equation (5.3). The formula (4.6) is applicable to the simple-layer potential $v_0(x)$, from which it follows that $v_0(x) = const$ in domain $R_m^+ \setminus \overline{\Omega}$. At infinity, a simple layer potential is zero, and consequently $v_0(x) \equiv 0$ in $R_m^+ \setminus \overline{\Omega}$, and also on the surface $\Gamma$. Applying now (4.5), we find that $v_0(x) \equiv 0$ is valid also inside the region $\Omega.$. But then $\frac{\partial}{\partial n}\{v_0(s)\}_i = 0,$ and by virtue of formula (4.4) we obtain $\rho_0(t) \equiv 0$. Thus, the homogeneous equation (5.3) has only a trivial solution; consequently, $\lambda = 2$ is not an eigenvalue of the kernel $K_1(s,t)$.

### 6. The Uniqueness Of Holmgren Problem's Solution

Let $\Omega$ be a finite domain in $R_m^+$ bounded by the open part $\Omega_1$ of the hyperplane $x_1 = 0$ and the Lyapunov surface $\Gamma$. We denote the boundary of the domain $\Omega_1$ by $\partial\Omega_1$.

**The Holmgren problem.** To find in domain $\Omega$ a regular solution of equation (1.1) that is continuous in the closed domain $\overline{\Omega}$ and satisfies the following boundary conditions

$$u\big|_\Gamma = \varphi(x) \quad (x \in \Gamma), \tag{6.1}$$

$$\lim_{x_1 \to 0} x_1^{2\alpha} \frac{\partial u}{\partial x_1} = \nu(x'), \quad (x' \in \Omega_1), \tag{6.2}$$

where $\varphi(x)$ and $\nu(x')$ are given continuous functions, and $x' = (x_2, x_3, ..., x_m)$.

**The uniqueness of the solution.** We consider the identity (2.1). Integrating both sides of this identity along the domain $\Omega_\varepsilon$ and using the Gauss-Ostrogradsky formula, we obtain

$$\int_{\Omega_\varepsilon} x_1^{2\alpha} \left[ u H_\alpha(v) - v H_\alpha(u) \right] dx = \int_{\partial\Omega_\varepsilon} x_1^{2\alpha} \sum_{i=1}^{m} \left( u \frac{\partial v}{\partial x_i} - v \frac{\partial u}{\partial x_i} \right) \cos(n, x_i) dS.$$

Here $\Omega_\varepsilon$ is a sub-domain of $\Omega$ at distance $\varepsilon$ from its boundary $\partial\Omega = \Omega_1 \cup \Gamma$. One can easily check that the following equality holds:

$$\int_{\Omega_\varepsilon} x_1^{2\alpha} u H_\alpha(u) dx = \int_{\Omega_\varepsilon} x_1^{2\alpha} \sum_{i=1}^{m} \left( \frac{\partial u}{\partial x_i} \right)^2 dx + \int_{\Omega_\varepsilon} \sum_{i=1}^{m} \frac{\partial}{\partial x_i} \left( x_1^{2\alpha} u \frac{\partial u}{\partial x_i} \right) dx.$$

An application of the Gauss-Ostrogradsky formula to this equality after $\varepsilon \to 0$ gives

$$\int_\Omega x_1^{2\alpha} \sum_{i=1}^{m} \left( \frac{\partial u}{\partial x_i} \right)^2 dx = \int_{\Omega_1} u(x) \nu(x) dS + \int_\Gamma x_1^{2\alpha} \varphi(x) \frac{\partial u}{\partial n} dS. \tag{6.3}$$

If we consider the homogeneous Holmgren problem, then from (6.3) one can get that

$$\int_\Omega x_1^{2\alpha} \sum_{i=1}^{m} \left( \frac{\partial u}{\partial x_i} \right)^2 dx = 0.$$

Hence, it follows that $u(x) = 0$ in $\overline{\Omega}$. Therefore the following uniqueness theorem holds:

**Theorem 5**. *If the Holmgren problem has a solution, then it is unique.*

### 7. Green's Function

**Definition.** We call $G(x, \xi)$ as Green's function of the Holmgren problem, if it satisfies the following conditions:

- this function is a regular solution of equation (1.1) in the domain $\Omega$, except at the point $\xi$, which is any fixed point of $\Omega$.
- it satisfies the boundary conditions

$$G(x, \xi)\big|_\Gamma = 0, \quad \frac{\partial G(x, \xi)}{\partial x_1}\bigg|_{x_1=0} = 0, \tag{7.1}$$

- it can be represented as
$$G(x, \xi) = q_1(x, \xi) + v_1(x, \xi), \tag{7.2}$$

where

$$q_1(x, \xi) = k_1 \left( r^2 \right)^{-\alpha - \frac{m-2}{2}} F\left( \alpha + \frac{m-2}{2}, \alpha; 2\alpha; \zeta \right)$$

is a fundamental solution of equation (1.1) and the function $v_1(x, \xi)$ regular solution of equation (1.1) in the domain $\Omega$.

The construction of the Green's function reduces to finding its regular part $v_1(x, x_0)$ which, by virtue of (7.1) and (7.2), must satisfy the boundary conditions

$$v_1(x, \xi)\big|_\Gamma = -q_1(x, \xi)\big|_\Gamma, \tag{7.3}$$

$$\left.\frac{\partial G(x,\xi)}{\partial x_1}\right|_{x_1=0}=0,$$

We look for the function $v_1(x,\xi)$ in the form of a double-layer potential:

$$v_1(x,\xi)=\int_\Gamma t_1^{2\alpha}\mu_1(t;\xi)\frac{\partial q_1(t,x)}{\partial n_t}d_t\Gamma. \tag{7.4}$$

Taking into account the equality (3.17) and the boundary condition (7.3), we obtain the integral equation for the denseness

$$\mu_1(s;\xi)-2\int_\Gamma \mu_1(t;\xi)K_1(s,t)d_t\Gamma=2q_1(s,\xi),\quad s\in\Gamma. \tag{7.5}$$

The right-hand side of (7.5) is a continuous function of $s$ (the point $\xi$ lies inside $\Omega$). In section 5 it was proved that it is not an eigenvalue of the kernel and, consequently, equation (7.5) is solvable and its continuous solution can be written in the form

$$\mu_1(s;\xi)=2q_1(s,\xi)+4\int_\Gamma R_1(s,t;2)q_1(t,\xi)d_t\Gamma, \tag{7.6}$$

where $R_1(s,t;\lambda)$ is a resolvent of kernel $K_1(s,t)$, $s\in\Gamma$. Substituting (7.6) into (7.4), we obtain

$$v_1(x,\xi)=2\int_\Gamma t_1^{2\alpha}q_1(t,\xi)\frac{\partial q_1(t,x)}{\partial n_t}d_t\Gamma$$
$$+4\int_\Gamma\int_\Gamma t_1^{2\alpha}R_1(t,s;2)q_1(s,\xi)\frac{\partial q_1(t,x)}{\partial n_t}d_t\Gamma d_s\Gamma. \tag{7.7}$$

We now define the function

$$g(x)=\begin{cases}v_1(x,\xi),\ x\in\Omega,\\ -q_1(x,\xi),\ x\in\Omega'.\end{cases} \tag{7.8}$$

The function $g(x)$ is a regular solution of (1.1) both inside the domain $\Omega$, and inside $\Omega'$ and equal to zero at infinity. Since point $\xi$ lies inside $\Omega$, then in $\Omega'$ the function $g(x)$ has derivatives of any order in all variables, continuous up to $\Gamma$. We can consider $g(x)$ in $\Omega'$ as a solution of equation (1.1) satisfying the boundary conditions

$$\left.\frac{\partial g(x)}{\partial n}\right|_\Gamma=-\left.\frac{\partial q_1(x;\xi)}{\partial n_x}\right|_\Gamma, \tag{7.9}$$

$$\left.\frac{\partial g(x)}{\partial x_1}\right|_{x_m=0}=0.$$

We represent this solution in the form of a simple-layer potential

$$g(x)=\int_\Gamma \rho_1(t;\xi)q_1(t,x)d_t\Gamma,\ x\in\Omega', \tag{7.10}$$

with an unknown denseness $\rho_1(t;\xi)$.

Using the formula (4.3), by virtue of condition (7.9), we obtain the integral equation for the denseness $\rho_1(t;\xi)$

$$\rho_1(s;\xi)-2\int_\Gamma \rho_1(t;\xi)K_1(t,s)d_t\Gamma=2\frac{\partial q_1(s,\xi)}{\partial n_s}. \tag{7.11}$$

Equation (7.11) is conjugated with equation (7.5). Its right-hand side is a continuous function of $s$. Thus, equation (7.11) has a continuous solution:

$$\rho_1(s;\xi)=2\frac{\partial q_1(s;\xi)}{\partial n_s}+4\int_\Gamma R_1(t,s;2)\frac{\partial q_1(t;\xi)}{\partial n_t}d_t\Gamma. \tag{7.12}$$

The values of a simple-layer potential $g(x)$ on the surface $\Gamma$ are equal to $-q_1(x;\xi)$, i.e. are the same as the functions $v_1(x;\xi)$, while on the hyperplane $x_1 = 0$ their partial derivatives with respect to $x_1$ are equal to zero.

From this, by the uniqueness theorem of the Holmgren problem, it follows that the formula (7.10) for the function $g(x)$ defined by (7.8) holds in the whole half-space $x_1 \geq 0$, that is,

$$v_1(x;\xi) = \int_\Gamma \rho_1(t;\xi) q_1(t,x) d_t\Gamma, \quad x \in \Omega. \tag{7.13}$$

Thus, the regular part $v_1(x;\xi)$ of the Green's function is representable in the form of a simple-layer potential.

Applying the formula (4.2) to (7.13), we obtain

$$2\frac{\partial}{\partial n_s}\{v_1(s;\xi)\}_i = \rho_1(s;\xi) + 2\int_\Gamma K_1(t,s)\rho_1(t;\xi) d_t\Gamma,$$

but, according to (7.11), we have

$$2\frac{\partial}{\partial n_s}\{q_1(s;\xi)\}_i = \rho_1(s;\xi) - 2\int_\Gamma K_1(t,s)\rho_1(t;\xi) d_t\Gamma.$$

Adding the last two equalities by term-wise and taking into account (7.2), we have

$$\frac{\partial}{\partial n_s}\{G_0(s;\xi)\} = 2\rho_1(s;\xi), \tag{7.14}$$

and, consequently, formula (7.13) can be written in the form

$$v_1(x;\xi) = \int_\Gamma q_1(t,x) \frac{\partial}{\partial n_t}\{G_1(t;\xi)\} d_t\Gamma. \tag{7.15}$$

Multiplying both sides of (7.12) by integrating over $s$ over the surface $\Gamma$ and by virtue of (7.6) and (7.4), we obtain

$$v_1(\xi;x) = \int_\Gamma \rho_1(t;\xi) q_1(t,x) d_t\Gamma.$$

Comparing this with formula (7.13), we have

$$v_1(\xi;x) = v_1(x;\xi), \tag{7.16}$$

if the points $x$ and $\xi$ are inside the domain $\Omega$.

**Lemma 3.** If points $x$ and $\xi$ are inside domain $\Omega$, then Green's function $G_1(x;\xi)$ is symmetric about those points.

The proof of the lemma 3 follows from the representation (7.2) of the Green's function and the equality (7.16).

For a normal domain $\Omega_0$ bounded by a hyperplane $x_1 = 0$ and a hemisphere

$$x_1^2 + x_2^2 + \ldots + x_m^2 = a^2,$$

the Green's function of the Holmgren problem has the form

$$G_{01}(x;\xi) = q_1(x;\xi) - \left(\frac{a}{R}\right)^{2\alpha_1} \cdot q_1\left(x; \frac{a^2}{R^2}\xi\right), \tag{7.17}$$

where

$$\xi_1^2 + \xi_2^2 + \ldots + \xi_m^2 = R^2.$$

We show that the function

$$v_{01}(x;\xi) = -\left(\frac{a}{R}\right)^{2\alpha_1} q_1(x;\xi)$$

can be represented in the form

$$v_{01}(x;\xi) = -\int_\Gamma \rho_1(s;x) v_{01}(s,\xi) d_s\Gamma,$$

where $\rho_1(s;x)$ is a solution of equation (7.11).

Indeed, let an arbitrary point $\xi$ be inside the domain $\Omega$. We consider the function
$$u(x;\xi) = -\int_\Gamma \rho_1(s;x) v_{01}(s,\xi) d_s\Gamma. \tag{7.18}$$
As a function of $x$, it satisfies (1.1), since this equation is satisfied by the function $\rho_1(s;x)$. Substituting expression (7.12) for $\rho_1(s;x)$, we obtain
$$u(x;\xi) = -\int_\Gamma \psi(s;\xi) \frac{\partial q_1(s;x)}{\partial n_s} d_s\Gamma, \tag{7.19}$$
where
$$\psi(s;\xi) = 2v_{01}(s,\xi) + 4\int_\Gamma R_1(s,t;2) v_{01}(t,\xi) d_t\Gamma,$$
i.e. $\psi(s;\xi)$ is a solution of equation
$$\psi(s;\xi) - 2\int_\Gamma K_1(s,t) \psi(t;\xi) d_t\Gamma = 2v_{01}(s;\xi). \tag{7.20}$$
Applying the formula (3.17) to the potential of the double layer (7.19), we obtain
$$u_i(s;\xi) = \frac{1}{2} \psi(s;\xi) - \int_\Gamma K_1(s,t) \psi(t;\xi) d_t\Gamma,$$
whence, by virtue of (7.20),
$$u_i(s;\xi) = v_{01}(s;\xi), \ \xi \in \Gamma.$$

It is easy to see that
$$\left.\frac{\partial u(x;\xi)}{\partial x_1}\right|_{x_1=0} = 0, \quad \left.\frac{\partial v_{01}(x;\xi)}{\partial x_1}\right|_{x_1=0} = 0.$$

Thus, the functions $u(x,\xi)$ and $v_{01}(x,\xi)$ satisfy the same equation (1.1) and the same boundary conditions, and by virtue of the uniqueness of the solution of the Holmgren problem, the equality $u(x,\xi) \equiv v_{01}(x,\xi)$ is satisfied.

Now, subtracting expression (7.17) from (7.2), we obtain
$$H_1(x,\xi) = G_1(x,\xi) - G_{01}(x,\xi) = v_1(x,\xi) - v_{01}(x,\xi)$$
or by virtue of (7.13), (7.16), (7.18) and (7.17) we obtain
$$H_1(x;\xi) = \int_\Gamma G_{01}(t;\xi) \rho_1(t;x) d_t\Gamma. \tag{7.21}$$

## 8. The solving the Holmgren problem for equation (1.1).

Let $x_0$ be a point inside the domain D. Consider the domain $\Omega_{\varepsilon,\delta} \subset \Omega$ bounded by the surface $\Gamma_\varepsilon$ which is parallel to the surface $\Gamma$, and the domain $\Gamma_{1\delta}$ lying on the hyperplane $x_1 = \delta > \varepsilon$. We choose $\varepsilon$ and $\delta$ so small that the point $\xi$ is inside $\Omega_{\varepsilon,\delta}$. We cut out from the domain $\Omega_{\varepsilon,\delta}$ a sphere of small radius $\rho$ with center at the point $x_0$ and the remainder part of $\Omega_{\varepsilon,\delta}$ denote by $\Omega_{\varepsilon,\delta}^\rho$, in which the Green's function $G_1(x;x_0)$ is a regular solution of (1.1).

Let $u(x)$ be a regular solution of equation (1.1) in the domain $\Omega$ that satisfies the boundary conditions (6.1) and (6.2). Applying the formula (1.2), we obtain
$$\int_{\Gamma_\varepsilon} x_1^{2\alpha}\left(G_1 \frac{\partial u}{\partial n} - u \frac{\partial G_1}{\partial n}\right) dS + \int_{\Gamma_{1\delta}} x_1^{2\alpha}\left(G_1 \frac{\partial u}{\partial n} - u \frac{\partial G_1}{\partial n}\right) dS = \int_{C_\rho} x_1^{2\alpha}\left(G_1 \frac{\partial u}{\partial n} - u \frac{\partial G_1}{\partial n}\right) dS,$$
Passing to the limit as $\rho \to 0$ and then as $\varepsilon \to 0$ and $\delta \to 0$, we obtain
$$u(x_0) = -\int_{\Gamma_1} \nu(x') G_1(0,x';x_0) dx' - \int_\Gamma \varphi(x) x_1^{2\alpha} \frac{\partial G_1(x,x_0)}{\partial n_x} d_x\Gamma, \tag{8.1}$$
We show that formula (8.1) gives a solution to the Holmgren problem.

It is easy to see that the first integral $I_1(x_0)$ in formula (8.1) is a solution of the equation (1.1) and is a regular in the domain $\Omega$, continuous in $\overline{\Omega}$.

We denote by

$$\vartheta(x_0) = \int_{\Gamma_1} v(x') q_1(0, x'; x_0) dx' = k_1 \int_{\Gamma_1} v(x') \left( x_{01}^2 + \sum_{i=2}^{m} (x_i - x_{0i})^2 \right)^{-\alpha - \frac{m-2}{2}} dx'. \qquad (8.2)$$

Here $\vartheta(x_0)$ is a continuous function in $\overline{\Omega}$. In view of (8.2) and (7.7) and the symmetry of the function $v_1(x; x_0)$, the integral $I_1(x_0)$ can be represented in the form

$$I_1(x_0) = -\vartheta(x_0) - 2\int_\Gamma \vartheta(t) t_1^{2\alpha} \frac{\partial G_1(t, x_0)}{\partial n_t} d_t\Gamma - 4\int_\Gamma \int_\Gamma R_1(t, s; 2)\vartheta(s) t_1^{2\alpha} \frac{\partial G_1(t, x_0)}{\partial n_t} d_t\Gamma d_s\Gamma. \qquad (8.3)$$

The last two integrals in the formula (8.3) are double-layer potentials. Taking into account the formula (3.17) and the integral equation for the resolvent $R_1(t, s; 2)$ from formula (8.3) we obtain

$$I_1(x_0)|_\Gamma = 0.$$

It is easy to see that

$$\lim_{x_{01} \to 0} x_{01}^{2\beta} \frac{\partial I_1(x_0)}{\partial x_0} = v(x_0'), \ x_0' \in \Omega_1.$$

In fact, the integral, by virtue of (7.13) and the symmetry of the function $v_1(x; x_0)$, can also be written in the form

$$\vartheta(x_0) = \int_{\Gamma_1} v(x') q_1(0, x'; x_0) dx' - \int_{\Gamma_1} v(x') dx' \int_\Gamma \rho(t; 0, x') q_1(t; x_0) dt.$$

It was shown in the section 6 that the derivative with respect to $x_{01}$ of the first term is equal to $v(x_0)$ at The derivative with respect to $x_{01}$ from the second term is zero at $x_{01} = 0$, since $\frac{\partial q_1}{\partial x_{01}} = 0$ at $x_{01} \to 0, \ x_0' \in \Omega_1$.

We consider the second integral $I_2(x_0)$ in the formula (8.1), which, by virtue of (7.14) and (7.12), can be written in the form

$$I_2(x_0) = -\int_\Gamma \varphi(t) \rho_1(t; x_0) d_t\Gamma = -\int_\Gamma \theta(t) t_1^{2\alpha} \frac{\partial}{\partial n_t} \{q_1(t, x_0)\} d_t\Gamma,$$

where

$$\theta(t) = 2\varphi(t) + 4\int_\Gamma R_1(t, s; 2)\varphi(s) d_s\Gamma,$$

i.e. the function $\theta(s)$ is a solution of the integral equation

$$\theta(s) - 2\int_\Gamma K_1(s; t)\theta(t) d_t\Gamma = 2\varphi(s). \qquad (8.4)$$

Since $\theta(s)$ is a continuous function, $I_2(x_0)$ is a solution of equation (1.1) regular in the domain $\Omega$ that is continuous in $\overline{\Omega}$, which, by virtue of (3.17), (3.18) and (8.4), satisfies condition

$$I_2(x_0)|_\Gamma = \varphi(s).$$

It is easy to see

$$\lim_{x_{01} \to 0} x_{01}^{2\alpha} \frac{\partial I_2(x_0)}{\partial x_0} = 0, \ x_0' \in \Omega_1.$$

**Comment.** Using the formulas (7.21) and (7.17), the solution (8.1) of the Holmgren problem for equation (1.1) can be written in the form

$$u(x_0) = -\int_{\Gamma_1} v(x') \left[ G_{01}(0, x'; x_0) + H_1(0, x'; x_0) \right] dx'$$

$$-\int_\Gamma \varphi(x) x_1^{2\alpha} \left[ \frac{\partial G_{01}(x, x_0)}{\partial n_x} + \frac{\partial H_1(x, x_0)}{\partial n_x} \right] d_x\Gamma, \qquad (8.5)$$

where

$$G_1(0, x'; x_0) = k_1 \left[ x_{01}^2 + \sum_{i=2}^{m}(x_i - x_{0i})^2 \right]^{-\alpha - \frac{m-2}{2}} - k_1 \left[ \sum_{i=2}^{m}\left(a - \frac{x_i x_{0i}}{a}\right)^2 + \frac{1}{a^2}\sum_{i=2}^{m} x_i^2 \sum_{j=1, j\neq i}^{m} x_{0j}^2 - (m-2)a^2 \right]^{-\alpha - \frac{m-2}{2}},$$

$$H_1(x; x_0) = \int_{\Gamma} \rho_1(t; x_0) G_{01}(t; x) d_t\Gamma.$$

The solution (8.5) of the Holmgren problem is more convenient for further investigations.

In the case of a hemisphere $\Omega_0$, the function $H_1(x; x_0) \equiv 0$ and the solution (8.5) take a simpler form:

$$u(x_0) = -k_1 \int_{\Gamma_1} \nu(x') \left\{ \left[ x_{01}^2 + \sum_{i=2}^{m}(x_i - x_{0i})^2 \right]^{-\alpha - \frac{m-2}{2}} \right.$$

$$\left. - \left[ \sum_{i=2}^{m}\left(a - \frac{x_i x_{0i}}{a}\right)^2 + \frac{1}{a^2}\sum_{i=2}^{m} x_i^2 \sum_{j=1, j\neq i}^{m} x_{0j}^2 - (m-2)a^2 \right]^{-\alpha - \frac{m-2}{2}} \right\} d_{x'}\Gamma_1 -$$

$$- (2\alpha + m - 2) k_1 (R^2 - a^2) \int_{\Gamma} \varphi(x) x_1^{2\alpha - 1} (r^2)^{-\alpha - \frac{m}{2}} F\left(\alpha + \frac{m}{2}, \alpha; 2\alpha; -\frac{4x_1 x_{01}}{r^2}\right) d_x\Gamma.$$

Thus, we obtain the solution of the Holmgren problem for the equation (1.1) in explicit form.